\documentclass[10pt]{article}
\usepackage{amsfonts}
\usepackage{amsmath}
\usepackage{mathrsfs}
\usepackage{mathrsfs,amscd,amssymb,amsthm,amsmath,bm,graphicx,psfrag,subfigure}

\makeatletter

\renewcommand{\@seccntformat}[1]{{\csname the#1\endcsname}{\normalsize .}\hspace{.5em}}
\makeatother

\def \[{\begin{equation}}
\def \]{\end{equation}}
\newtheorem{thm}{Theorem}[section]

\newtheorem{defi}{Definition}

\newtheorem{fac}{Fact}

\newtheorem{lem}[thm]{Lemma}

\newenvironment{wst}
{\setlength{\leftmargini}{1.5\parindent}
 \begin{itemize}
 \setlength{\itemsep}{-1.1mm}}
{\end{itemize}}

\begin{document}
\setlength{\baselineskip}{15pt}
\begin{center}{\Large \bf
Extremal values on the eccentric distance sum of trees\footnote{Financially supported by the National Natural Science
Foundation of China (Grant No. 11071096)  and the Special Fund for Basic Scientific Research of Central Colleges (CCNU11A02015).}}

\vspace{2mm}

{\large Shuchao Li\footnote{E-mail: lscmath@mail.ccnu.edu.cn (S.C.
Li), zizaifei001@163.com (M. Zhang)}, Meng Zhang}\vspace{2mm}

Faculty of Mathematics and Statistics,  Central China Normal
University, Wuhan 430079, P.R. China\vspace{1mm}
\end{center}

\noindent {\bf Abstract}:\ Let $G=(V_G, E_G)$ be a simple connected
graph. The eccentric distance sum of $G$ is defined as $\xi^{d}(G) =
\sum _{v\in V_G}\varepsilon_{G}(v)D_{G}(v)$, where
$\varepsilon_G(v)$ is the eccentricity of the vertex $v$ and $D_G(v)
= \sum_{u\in V_G}d_G(u,v)$ is the sum of all distances from the
vertex $v$. In this paper the tree
among $n$-vertex trees with domination number $\gamma$ having the
minimal eccentric distance sum is determined and the tree
among $n$-vertex trees with domination number $\gamma$ satisfying
$n = k\gamma$ having the maximal eccentric distance sum
is identified, respectively, for
$k=2,3,\frac{n}{3},\frac{n}{2}$. Sharp upper and lower bounds on the
eccentric distance sums among the $n$-vertex
trees with $k$ leaves are determined. Finally, the trees
among the $n$-vertex trees with a given
bipartition having the minimal, second minimal and third minimal
eccentric distance sums are determined, respectively.

\vspace{2mm} \noindent{\it Keywords}: Eccentric distance sum;
Domination number; Leaves; Bipartition

\vspace{2mm}


{\setcounter{section}{0}
\section{\normalsize Introduction}\setcounter{equation}{0}

We consider only simple connected graphs.  Let $G=(V_G, E_G)$ be a graph
with $u,v\in V_G$, $d_G(u)$ (or $d(u)$ for short) denotes the
degree of $u$; we call $u$ a \textit{leaf} if $d_G(u)=1$. Let $PV(G)$ be the set of all leaves of $G$. The \textit{distance} $d_G(u,v)$ is defined as the
length of the shortest path between $u$ and $v$ in $G$; $D_G(u)$ (or
$D(u)$ for short) denotes the sum of distances between $u$ and all other
vertices of $G$. Let
$N_G(u)$ be the set of vertices adjacent to $u$ in $G$.  The \textit{eccentricity}
$\varepsilon(v)$ of a vertex $v$ is the maximum distance from $v$ to
any other vertex. The \textit{radius} $rad(G)$ of a graph is the
minimum eccentricity of any vertex, while the diameter $diam(G)$ of
a graph is the maximum eccentricity of any vertex in the graph. The
\textit{center} of a graph is the vertices whose eccentricity is
equal to the radius. $G-u$ denotes the graph obtained from $G$ by
deleting the vertex $u$ together with its incident edges (this
notation is naturally extended if more than one vertex is deleted).
If $U\subseteq V_G$, then $G[U]$ denotes the graph on $U$ whose edges are precisely the edges of $G$ with both ends in $U$.
Let $S_n$, $P_n$ and $K_n$ be a star, a path and a complete graph on
$n$ vertices, respectively. For a real number $x$ we denote by $\lfloor x\rfloor$
the greatest integer no greater than $x$, and by $\lceil x\rceil$ the least integer no less than $x$.

A single number that can be used to characterize some property of
the graph of a molecule is called a \textit{topological index}, or
\textit{graph invariant}. Topological index is a graph theoretic
property that is preserved by isomorphism. The chemical information
derived through topological index has been found useful in chemical
documentation, isomer discrimination, structure property
correlations, etc. \cite{0}. For quite some time there has been
rising interest in the field of computational chemistry in
topological indices. The interest in topological indices is mainly
related to their use in nonempirical quantitative structure-property
relationships and quantitative structure-activity relationships.
Among various indices, the \textit{Wiener index} has been one of the
most widely used descriptors in quantitative structure activity
relationships. Many recently established topological indices such as
\textit{degree distance index}, \textit{eccentric connectivity
index} and so on are used as molecular descriptors.

The \textit{Wiener index} is defined as the sum of all distances
between unordered pairs of vertices
$$
  W(G) =\sum_{u,v\in V_G}d_G(u, v).
$$
It is considered as one of the most used topological index with high
correlation with many physical and chemical properties of a molecule
(modelled by a graph). For the recent survey on Wiener index one may
refer to \cite{4} and the references cited in.

The  \textit{degree distance index} $DD(G)$  was introduced by
Dobrynin and Kochetova \cite{05} and Gutman \cite{09} as
graph-theoretical descriptor for characterizing alkanes; it can be
considered as a weighted version of the Wiener index
$$
DD(G) = \sum_{u,v\in V_G}(d_G(u)+ d_G(v))d_G(u, v) = \sum_{v\in
V_G}d_G(v)D_G(v),
$$
where the summation goes over all pairs of vertices in $G$. 

Sharma, Goswami and Madan \cite{18} introduced a distance-based
molecular structure descriptor, \textit{eccentric connectivity
index} (ECI) defined as
$$
\zeta^c(G) =\sum_{v\in V_G} \varepsilon_G(v)d_G(v).
$$
The index $\zeta^c(G)$ was successfully used for mathematical models
of biological activities of diverse nature \cite{6,7}. For the study
of its mathematical properties one may be referred to \cite{A1,A0,15} and
the references there in.

Recently, a novel graph invariant for predicting biological and
physical properties-\textit{eccentric distance sum} (EDS) was
introduced by Gupta, Singh and Madan \cite{S0}, which was defined as
$$
\xi^d(G) =\sum_{v\in V_G} \varepsilon_G(v)D_G(v).
$$
The eccentric distance sum can be defined alternatively as
$$
\xi^d(G) =\sum_{u,v\in V_G} (\varepsilon_G(v)+\varepsilon_G(u))d_G(u,v).
$$
This topological index has vast potential in structure
activity/property relationships; it also displays high
discriminating power with respect to both biological activity and
physical properties; see \cite{S0}. From \cite{S0} we also know
that some structure activity and quantitative structure property
studies using eccentric distance sum were better than the
corresponding values obtained using the Wiener index. It is also
interesting to study the mathematical property of this topological
index. Yu, Feng and Ili\'{c} \cite{12} identified the extremal
unicyclic graphs of given girth having the minimal and second
minimal EDS; they also characterized the trees with the minimal EDS
among the $n$-vertex trees of a given diameter. Hua, Xu and Shu
\cite{H-H-X-K} obtained the sharp lower bound on EDS of $n$-vertex
cacti. Hua, Zhang and Xu \cite{3} studied the graphs with graph parameters having the minimum EDS.
Ili\'{c}, Yu and Feng  \cite{13} studied the various lower
and upper bounds for the EDS in terms of the other graph invariant
including the Wiener index, the degree distance index, the eccentric
connectivity index and so on. Yu, Feng and authors here \cite{17} identified the trees
with the minimal and second minimal eccentric distance sums among
the $n$-vertex trees with matching number $q$; as well they characterized the extremal tree with
the second minimal eccentric distance sum among the $n$-vertex trees of a given diameter.
Consequently, they determined the trees with the third and fourth minimal eccentric distance
sums among the $n$-vertex trees. Motivated by
these results it is natural for us to continue
the research on the eccentric distance sum of trees.

This paper is organized as follows. We first characterize the
trees with the minimal EDS among $n$-vertex trees with
domination number $\gamma$, as well we determine the trees with the
maximal EDS among $n$-vertex trees with domination number $\gamma$
satisfying $n = k\gamma$, where $k = 2,3,\frac{n}{3},\frac{n}{2}$. Then we
identify the trees with the minimal and maximal EDS among the
$n$-vertex trees with $k$ leaves, respectively. Finally, we characterize
trees with the minimal, second minimal and third minimal EDS among the $n$-vertex trees of a given bipartition
$(p,q)$.

\section{\normalsize The extremal EDS of vertex trees with domination number $\gamma$}
In this section, we characterize the tree with
the minimal EDS among $n$-vertex trees with domination number
$\gamma$; as well we determine the tree with maximal EDS among $n$-vertex
trees with domination number $\gamma$ for $n = k\gamma,\, k =
2,3,\frac{n}{3},\frac{n}{2}$.  For convenience, let $\mathscr{T}_{n,\gamma}$ be
the set of all $n$-vertex trees with domination number $\gamma.$

\begin{lem}[\cite{1}]
For a graph $G$, we have $\gamma(G) \leq \beta(G)$.
\end{lem}

In \cite{8} and in many subsequent works (see especially
\cite{10,11}) it has been demonstrated the following conclusion.

\begin{lem}
Let $T$ be a tree of order $n$. Then
$W(S_n) \leq W(T) \leq W
(P_n) .$
The left equality holds if and only if $T\cong S_n$, and the right
equality holds if and only if $T\cong P_n$.
\end{lem}

\begin{lem}[\cite{13,12}]
Let $T$ be a tree of order $n$. Then
$\xi^d(S_n) \leq \xi^d(T) \leq \xi^d(P_n) .$
The left equality holds if and only if $T\cong S_n$, and the right
equality holds if and only if $T\cong P_n$.
\end{lem}

Let $P_l(a,b)$ be an $n$-vertex tree obtained by attaching $a$ and $b$
leaves to the two endvertices of $P_l=v_1v_2\ldots
v_l,(l\geq2)$, respectively. Here, $a+b=n-l,\ a,b\geq 1.$
\begin{lem}
$
\xi^d(P_l(1,n-l-1))<\xi^d(P_l(2,n-l-2))<\cdots<\xi^d(P_l(\lfloor\frac{n-l}{2}\rfloor,\lceil\frac{n-l}{2}\rceil)).
$
\end{lem}
\begin{proof}
It suffices to show that $\xi^d(P_l(a-1,b+1))>\xi^d(P_l(a,b))$ if
$a-b>1$.

Let $\varepsilon'(x)$ (resp. $\varepsilon(x)$) denote the
eccentricity of $x$ in $P_l(a-1,b+1)$ (resp. $P_l(a,b)$),\, $D'(x)$
(resp. $D(x)$) denote the sum of all distances from the vertex $x$ in $P_l(a-1,b+1)$ (resp.
$P_l(a,b)$). It is obvious to see that the eccentricity of every
vertex remains the same. By the definition of EDS, we get
\begin{eqnarray*}
  \xi^d(P_l(a-1,b+1))-\xi^d(P_l(a,b)) &=& \sum_{x\in V_{P_l(a,b)}}\varepsilon(x)(D'(x)-D(x)) \\
   &=& \sum_{i=1}^{l}\varepsilon(v_i)(l+1-2i)+(a-1)(\varepsilon(v_1)+1)(l-1)
   +b(\varepsilon(v_l)+1)(1-l)\\
   &&+(l+1)(l-1)(a-b-1) \\
   &=& \sum_{i=1}^{l}\varepsilon(v_i)(l+1-2i)+2(l+1)(l-1)(a-b-1).
\end{eqnarray*}
Note that $l-1>0$, $a-b-1\geq1$ and
$\sum_{i=1}^{l}\varepsilon(v_i)(l+1-2i)\geq0$. Hence,
$\xi^d(P_l(a-1,b+1))>\xi^d(P_l(a,b))$, as desired.
\end{proof}

Let $T$ be a tree of order $n>3$ and $e=uv$ be a nonpendant edge.
Suppose that $T-e=T_1\cup T_2$ with $u\in V_{T_1}$ and $v\in V_{T_2}$.
Now we construct a new tree $T_0$ obtained by identifying the vertex
$u$ of $T_1$ with vertex $v$ of $T_2$ and attaching a leaf
to the $u(= v)$. Then we say that $T_0$ is obtained by running
\textit{edge-growing transformation} of $T$ (on edge $e=uv$), or
e.g.t of $T$ (on edge $e =uv$) for short; see Fig. 1.
\begin{figure}[h!]
\begin{center}
  \psfrag{a}{$u$}\psfrag{b}{$v$}\psfrag{c}{$u(=v)$}\psfrag{d}{e.g.t}\psfrag{e}{$T$}
  \psfrag{f}{$T_0$}\psfrag{1}{$T_1$}\psfrag{2}{$T_2$}
  \includegraphics[width=120mm] {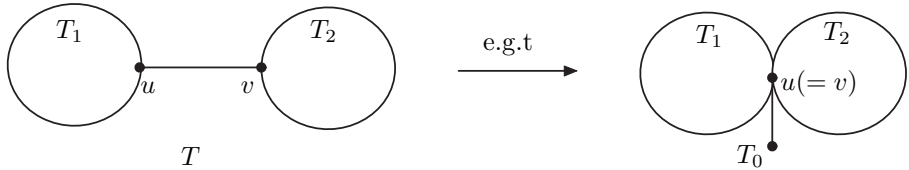}\\
  \caption{Tow trees $T$ and $T_0$ for e.g.t.}
\end{center}
\end{figure}

\begin{lem}[\cite{3}]
Let $T$ be a tree of order $n>3$ and $e=uv$ be a nonpendant edge of
$T$. If $T_0$ is a tree obtained from $T$ by running one step of
e.g.t (on edge $e=uv$), then we have $\xi^d(T_0)< \xi^d(T)$.
\end{lem}

As introduced in \cite{14}, we use $\xi(G)=\sum_{x\in
V_G}\varepsilon_G(x)$ to denote the \textit{total eccentricity} of
a connected graph $G$. The corona of two graphs $G_1$ and $G_2$, introduced in \cite{16},
is a new graph $G=G_1\circ G_2$ obtained from one copy of $G_1$ with
$|V_{G_1}|$ copies of $G_2$ where the $i$th vertex of $G_1$ is
adjacent to every vertex in the $i$th copy of $G_2$. As an example,
the corona $G\circ K_1$ is a graph obtained from attaching a leaf to each vertex of $G$. In particular, for a positive integer
$p$, we denote by $G^{(p)}$ the graph obtained by attaching $p$
leaves to every vertex of $G$. Note that $G^{(p)}$ has $(p
+ 1)n$ vertices and $G\circ K_1=G^{(1)}$.

\begin{lem}
Let $T$ be a tree of order $n$ and $T^{(m)}$ be the graph as defined
above. Then
$\xi^d(S_n^{(m)}) \leq \xi^d(T^{(m)}) \leq \xi^d(P_n^{(m)}) .$
The left equality holds if and only if $T\cong S_n$, and the right
equality holds if and only if $T\cong P_n$.
\end{lem}
\begin{proof}
Let $uv$ be a pendant edge of $T^{(m)}$ with $d_{T^{(m)}}(u) =1$,
then $\varepsilon_{T^{(m)}}(u) = \varepsilon_{T^{(m)}}(v) + 1$
and $D_{T^{(m)}}(u) = D_{T^{(m)}}(v) + (m+1)n - 2$.
Moreover,
$$
\text{$\varepsilon_{T^{(m)}}(x) = \varepsilon_T(x) + 1$\ \ and\ \
$D_{T^{(m)}}(x) = (m+1)D_{T}(x) + nm$ \ \ for any vertex $x\in V_T$.}
$$
By the definition of EDS, we have
\begin{align*}
  \xi^d(T^{(m)}) =& \sum_{x\in V_{T^{(m)}}}\varepsilon _{T^{(m)}}(x)D_{T^{(m)}}(x) \\
    =& \sum_{x\in V_T}\varepsilon _{T^{(m)}}(x)D_{T^{(m)}}(x) + m\sum_{x\in V_T}(\varepsilon_{T^{(m)}}(x) + 1)(D_{T^{(m)}}(x)+(m+1)n-2) \\
    =&\sum_{x\in V_T}(\varepsilon _T(x)+1)((m+1)D_{T}(x) + nm) + m\sum_{x\in V_T}(\varepsilon_T(x) + 2)((m+1)D_{T}(x) +2nm+n-2)   \\
    =& (m+1)^2\sum_{x\in V_T}\varepsilon_T(x)D_T(x)+(2m+1)(m+1)\sum_{x\in V_T}D_T(x)+2m(nm+n-1)\sum_{x\in V_T}\varepsilon_T(x)\\
     &+4n^2m^2+3n^2m-4nm \\
    =&(m+1)^2\xi^d(T)+2(2m+1)(m+1)W(T)+2m(nm+n-1)\xi(T)+4n^2m^2+3n^2m-4nm.
\end{align*}

By Lemmas $2.2$ and $2.3$, it suffices to show that
$\xi(S_n) \leq \xi(T) \leq \xi(P_n),$
the left equality holds if and only if $T\cong S_n$, whereas the
right equality holds if and only if $T\cong P_n$.

In fact, assume to the contrary that $T\ncong S_n$, then we have
$\varepsilon_T(x)\geq 2$ for every vertex $x\in
V_T$; otherwise, $T\cong S_n$. Hence,  $\xi(T)\geq
2n>2n-1=\xi(S_n)$

On the other hand, suppose that $P$ is one of the longest paths of $T$, and let
$d={\rm diam}(T)=|V_P|-1.$  Note that for every vertex $x\in V_T$,
$\varepsilon_T(x)\leq d$, hence
\begin{eqnarray*}
     &&\xi(T) = \sum_{x\in V_T}\varepsilon_T(x)
       \leq \xi(P)+(n-d-1)d
       =\left\{
            \begin{array}{ll}
              nd-\frac{1}{4}d^2=f(d), & \hbox{$d$ is even;} \\[5pt]
              nd-\frac{1}{4}d^2+\frac{1}{4}=g(d), & \hbox{$d$ is odd.}
            \end{array}
          \right.
   \end{eqnarray*}
It is routine to check that $f(d)$ (resp. $g(d)$) is a strictly increasing
function in $d$ with $d\in [2,n-1]$.

Note that
\begin{eqnarray*}
  f(d) &=& nd-\frac{1}{4}d^2<n(d+1)-\frac{1}{4}(d+1)^2+\frac{1}{4}=g(d+1), \\
  g(d) &=& nd-\frac{1}{4}d^2+\frac{1}{4}<n(d+1)-\frac{1}{4}(d+1)^2=f(d),
\end{eqnarray*}
we have $\xi(T)\leq\xi(P_n)$, and the equality holds if and only if
$d=n-1$, i.e., $T\cong P_n$.

This completes the proof.
\end{proof}
A subset $S$ of $V_G$ is called a \textit{dominating set} of $G$ if
 for every vertex $v\in V_G\setminus S$, there exists a vertex $u\in S$ such that
$v$ is adjacent to $u$. A vertex in the dominating set is called
\textit{dominating vertex}. For a dominating set $S$ of graph $G$
with $v\in S$, $u\in V_G\setminus S$, if $vu\in E_G$, then $u$ is said to
be dominated by $v$. The \textit{domination number} of $G$, denoted
by $\gamma(G)$, is defined as the minimum cardinality of dominating
sets of $G$. For a connected graph $G$ of order $n$, Ore \cite{19}
obtained that $\gamma(G)\leq\frac{n}{2}$. And the equality case was
characterized independently in \cite{20,21}. For a graph $G$, the
\textit{matching number} $\beta(G)$ is the cardinality of a maximum
matching of $G$. We denote by $T_{n,\beta}$ the tree obtained from the star
$S_{n-\beta+1}$ by attaching a pendant edge to each of certain
$\beta-1$ non-central vertices of $S_{n-\beta+1}$. It is easy to see
that $T_{n,\beta}$ contains an $\beta$-matching. If $n=2\beta$, then
it has a perfect matching. Tree $T_{n,\beta}$ is depicted in Fig. 2.
\begin{figure}[h!]
\begin{center}
  \psfrag{a}{$p-1$}\psfrag{b}{$q-1$}\psfrag{c}{$T(p,q)$}\psfrag{d}{$n-2\beta+1$}\psfrag{e}{$\beta-1$}\psfrag{f}{$T_{n,\beta}$}
  \includegraphics[width=60mm]{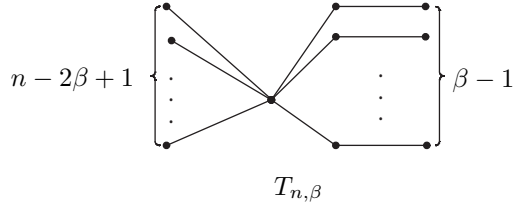}\\
  \caption{Trees $T_{n,\beta}$.}
\end{center}
\end{figure}

\begin{lem}[\cite{17}]
Among all the trees of order $n$ and with the matching number
$\beta$, the tree $T_{n,\beta}$ has the minimal EDS, and
$\xi^d(T_{n,\beta})=6n^{2}+\beta^{2}+9\beta n-22n-28\beta+34$.
\end{lem}

\begin{lem}[\cite{20,21}]
If $n=2\gamma$, then a tree $T$ belongs to $\mathscr{T}_{n,\gamma}$
if and only if there exists a tree $H$ of order $\gamma$ such that
$T=H\circ K_1$.
\end{lem}

\begin{lem}
If $T'\in \mathscr{T}_{n,\gamma}$ has the minimal EDS, then we have
$\gamma(T')=\beta(T')=\gamma$.
\end{lem}
\begin{proof}
By Lemma $2.1$, it suffices to show that
$\gamma(T')\geq\beta(T')$.
\begin{figure}[h!]
\begin{center}
  \psfrag{1}{$v_1$}\psfrag{2}{$v_{1'}$}\psfrag{3}{$v_2$}\psfrag{4}{$v_{2'}$}\psfrag{5}{$v_\gamma$}
  \psfrag{6}{$v_{\gamma'}$}\psfrag{7}{$w_1$}\psfrag{8}{$w_2$}\psfrag{a}{$T'$}\psfrag{b}{$T''$}
  \includegraphics[width=100mm] {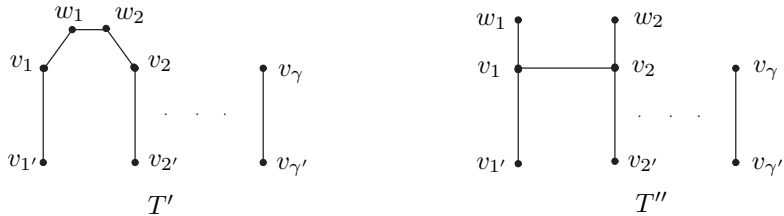}\\
  \caption{The structures of $T'$ and $T''$ }
\end{center}
\end{figure}
Otherwise, by the definition of the set $\mathscr{T}_{n,\gamma}$, we
have $\gamma=\gamma(T')<\beta(T')$. Assume that
$S=\{v_1,v_2,\cdots,v_\gamma\}$ is a dominating set of cardinality
$\gamma$. Then there exist $\gamma$ independent edges
$v_1v'_1,v_2v'_2,\ldots,v_\gamma v'_\gamma$ in $T_0$. Note that
$\gamma=\gamma(T')<\beta(T')$, there must exist another edge, say
$w_1w_2$, which is independent of each edge $v_iv'_i,$
$i=1,2,\ldots,\gamma$.

If $w_1,\, w_2$ are dominated by the same vertex $v_i\in
S$, then a triangle $C_3=w_1w_2v_i$ occurs. This is impossible
because of the fact that $T'$ is a tree. Therefore we claim that
two vertices $w_1,w_2$ are dominated by two deferent vertices from $S$.
Without loss of generality, assume that $w_i$ is dominated by the
vertex $v_i$ for $i=1, 2$ (see Fig.~3). Now we construct a new tree
$T''\in \mathscr{T}_{n,\gamma}$ by running e.g.t. of $T_0$ on the
edges $v_1w_1$ and $v_2w_2$, respectively. By Lemma $2.5$, we have
$\xi^d(T'')<\xi^d(T')$. This contradicts the choice of $T'$.
Thus we complete the proof of this lemma.
\end{proof}

Combining Lemmas 2.7 and 2.9, the following is obvious.
\begin{thm}
For any tree $T\in \mathscr{T}_{n,\gamma}$, we have
$
\xi^d(T)\geq 6n^{2}+\gamma^{2}+9\gamma n-22n-28\gamma+34.
$
The equality holds if and only if $T\cong T_{n,\gamma}$.
\end{thm}

\begin{thm}
Among all the trees from $\mathscr{T}_{n,\frac{n}{2}},$ the
tree $P_{\frac{n}{2}}\circ K_1$ has the maximal EDS.
\end{thm}
\begin{proof}
By Lemma 2.8, any tree from $\mathscr{T}_{n,\frac{n}{2}}$ must be of
the form $H\circ K_1$ where $H$ is a tree of order
$\frac{n}{2}=\gamma$. Taking $m=1$ in Lemma 2.6 implies our result immediately.
\end{proof}

\begin{thm}
Among all the trees in $\mathscr{T}_{n,\lceil\frac{n}{3}\rceil}$
with $n>4$, the tree $P_n$ has the maximal EDS.
\end{thm}

\begin{proof}
It is known \cite{13} that the path  $P_n=v_1v_2\ldots v_n$ has
the maximal EDS among all the trees of order $n$. Hence, in order to complete the proof, it suffices to show that
$\gamma(P_n) =\lceil \frac{n}{3}\rceil $.

Assume that $n=3k+r,\, 0 \leq r \leq 2$ and let $S_0 = \{v_2, v_5, \ldots , v_{3k-1}\}$.
Note that the vertex subset $S_0$ (resp. $S_0 \cup\{v_{3k+1}\}$) is a dominating set of $P_n$ for
$n = 3k$ (resp. $n=3k+1$,\, $3k + 2$). By
the definition of the domination number, we have $\gamma(P_n) \leq\lceil
\frac{n}{3}\rceil $. If $\gamma(P_n)<\lceil \frac{n}{3}\rceil $ ,
that is, $\gamma(P_n) \leq\lceil \frac{n}{3}\rceil -1$, then we
claim that at least three vertices are dominated by one vertex from a
dominating set. By the structure of $P_n$, this is impossible. So we
have $\gamma(P_n) =\lceil \frac{n}{3}\rceil$, as desired. 
\end{proof}

\begin{thm}
Among all the trees from $\mathscr{T}_{n,2}$ with $n \geq4$, the
tree $P_4(\lfloor\frac{n-4}{2}\rfloor,\lceil\frac{n-4}{2}\rceil)$ has
the maximal EDS.
\end{thm}
\begin{proof}
In view of Theorem 2.12, our result holds for $n = 4, 5, 6$. So in what follows
we only consider the case for $n\geq 7$. Assume that $T_1
\in\mathscr{T}_{n,2}$ has the maximal EDS and $S =\{w_1, w_2\}$
is a dominating set of $T_1$. Now we show the following two
claims:\vspace{2mm}

\noindent\textbf{Claim 1.}\  $w_1$ is not adjacent to $w_2$ .\vspace{2mm}

\noindent{\bf Proof of Claim 1.}\ \
If not, then $T_1$ must be of the form $P_2(a,b)$ with $a+b=n-2$ and $a \leq b$. By Lemma 2.4,
we have $b-a \leq 1$. That is to say, $T_1 \cong
P_2(\lfloor\frac{n-2}{2}\rfloor,\lceil\frac{n-2}{2}\rceil)$. Note
that $\frac{n-2}{2} \geq \frac{5}{2} > 2$. After running the
converse of e.g.t. on the edge $w_1w_2$ of $T_1$, we obtain a new
tree $T_2 \cong
P_3(\lfloor\frac{n-2}{2}\rfloor,\lceil\frac{n-2}{2}\rceil-1)$ which
still belongs to $\mathscr{T}_{n,2}$. By Lemma 2.5, we have
$\xi^d(T_2)> \xi^d(T_1)$, which contradicts the choice of $T_1$.\qed\vspace{2mm}

\noindent\textbf{Claim 2.}\  $d_{T_1}(w_1,w_2)=3$.\vspace{2mm}

\noindent{\bf Proof of Claim 2.}\ \ From Claim 1, we have $d_{T_1}(w_1,w_2)\geq 2$. If $d_{T_1}(w_1,w_2)\geq 4$, then
there exists at least one vertex $x$ on the shortest path between
$w_1$ and $w_2$ such that $x$ can not be dominated by the two
vertices $w_1$ and $w_2$. This contradicts the fact that
$T_1\in\mathscr{T}_{n,2}$. Then we get $2\leq d_{T_1}(w_1,w_2)\leq 3$. If
$d_{T_1}(w_1,w_2)= 2$, then we find that $T_1 \cong
P_3(\lfloor\frac{n-3}{2}\rfloor,\lceil\frac{n-3}{2}\rceil)$ by Lemma 2.4. Assume
that the common neighbor of $w_1$ and $w_2$ is $w_0$. Note that
$\frac{n-3}{2}\geq 2> 1$. By running the converse of e.g.t. on the
edge $w_0w_1$ or $w_0w_2$, in view of Lemma 2.5, we get a new tree
of the form $P_4(a, b)$ with $a+b =n-4$, which is still in
$\mathscr{T}_{n,2}$ but has a larger EDS. This is impossible because
of the maximality of $\xi^d(T_1)$, as desired. \qed

By Claims 1 and 2, $T_1$ must be of the form $P_4(a,b)$
with $a+b=n-4$. By Lemma 2.4, this theorem follows immediately.
\end{proof}

\section{\normalsize The extremal EDS of trees with $k$ leaves}
\setcounter{equation}{0}

In this section, we are to determine the trees with the minimal and
maximal EDS among the $n$-vertex trees each of which contains $k$ leaves. Note that
there is just one tree for $k=n-1$ or 2, hence in what follows we
consider $3\le k\leq n-2.$ For convenience, let $\mathscr{T}_n^k$ be
the set of all $n$-vertex trees with $k$ leaves.

A \textit{spider} is a tree with at most one vertex of degree more than 2, called the \textit{hub} of the spider (if no vertex of degree more than two, then any vertex can be the hub). A \textit{leg} of a spider is a path from the hub to one of its leaves. Let $S(a_1,a_2,\ldots,a_k)$ be a spider with $k$ legs $P^1,P^2,\ldots,P^k$ satisfying the length of $P^i$ is $a_i$
$(i=1,2,\ldots,k)$, and $\sum_{i=1}^ka_i=n-1$. Call $S(a_1,a_2,\ldots,a_k)$ a \textit{balanced
spider} if $\mid a_i-a_j\mid\leq 1$ for $1 \leqslant i, j \leqslant k.$
\begin{defi}
Let $T$ be an arbitrary tree rooted at a center vertex and let $v$
be a vertex of degree $m+1$ ($m\geq2$). Suppose that $w$ is
adjacent to $v$ with $\varepsilon_T(v)\geq\varepsilon_T(w)$ and
that $T_1,T_2,\ldots,T_m$ are subtrees under $v$ with root vertices
$v_1,v_2,\ldots,v_m$ such that the tree $T_m$ is actually a path. Let
$T'=T-\{vv_1,vv_2,\ldots,vv_{m-1}\}+\{wv_1,wv_2,\ldots,wv_{m-1}\}$. We
say that $T'$ is a $\rho$ transformation of $T$ and denote it by $T'=\rho(T,v)$ (see Fig. 4).
\end{defi}

Note that, by Definition 1, $|PV(T')|=|PV(T)|$ and ${\rm diam}(T')\le {\rm diam}(T).$
\begin{figure}[h!]
\begin{center}
  \psfrag{s}{$G_0$}
  \psfrag{y}{$w$}\psfrag{x}{$u$}\psfrag{z}{$v$}\psfrag{1}{$v_1$}\psfrag{2}{$v_2$}\psfrag{3}{$v_{m-1}$}\psfrag{4}{$v_m$}
  \psfrag{a}{$T_1$}\psfrag{b}{$T_2$}\psfrag{c}{$T_{m-1}$}\psfrag{d}{$T_m$}\psfrag{8}{$T$}\psfrag{9}{$T'$}\psfrag{e}{\vdots}
  \includegraphics[width=150mm] {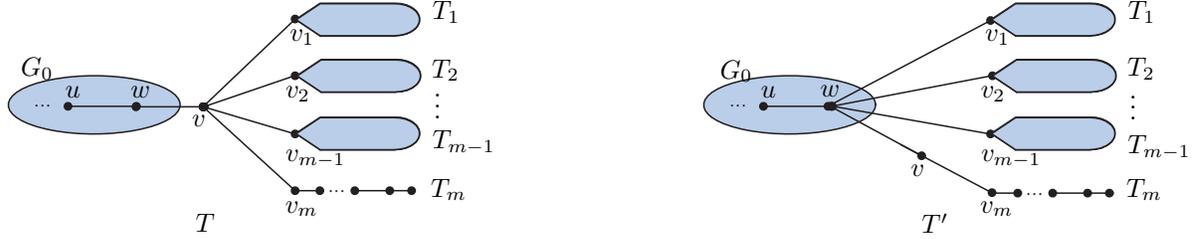}\\
  \caption{$\rho$ Transformation }
\end{center}
\end{figure}

\begin{lem}
Let $T$ and $T'$ be the trees defined as above, one has
$
\xi^d(T)\geq \xi^d(T').
$
The equality holds if and only if
$\varepsilon_T(v)=\varepsilon_T(w)$ and $T[S]$ is one of the longest paths in $T$, where
$S=V_{G_0}\cup V_{T_m}\cup \{v\}.$
\end{lem}
\begin{proof}
Let $G_0$ be the graph obtained from $T$ (or $T'$) by deleting $V_{T_1}\cup V_{T_2} \cup \ldots\cup V_{T_m}\cup \{v\}$ (see Fig. 4).
Let $P^i$ be one of the longest paths contained in $T_i$ such that one of its endvertices is $v_i, \, i=1,2,\ldots, m-1.$
From the definition of $\rho$ transformation, $w$ is adjacent to $v$
with $\varepsilon_T(v)\geq\varepsilon_T(w)$, we have the following fact.
\begin{fac}
$G_0$ contains one center, say $c$, of $T$ and there exists a longest path $P=wu_1u_2\ldots c \ldots y$ in $G_0$ such that $y\in PV(T)$ and $|V_P| >|V_{T_m}|$.
\end{fac}

Moreover, we have
\begin{eqnarray}
 \varepsilon_T(x)&\geq&\varepsilon_{T'}(x)>0\ \ \ \text{for all $x\in V_{G_0}$,}\label{eq:0031} \\
  \varepsilon_T(x)&=&\varepsilon_{T'}(x)>0\ \ \ \text{for all $x\in V_{T_m}\cup \{v\}$,}\label{eq:0032} \\
\varepsilon_T(x)&=&\varepsilon_{T}(v)+d_T(x,v)=\varepsilon_{T'}(v)+d_{T'}(x,w)\notag\\
         &\geq&\varepsilon_{T'}(w)+d_{T'}(x,w)=\varepsilon_{T'}(x)\ \ \ \ \ \ \text{for all $x\in V_{T_1}\cup V_{T_2}\cup \ldots\cup V_{T_{m-1}}$,}\label{eq:0033}\\
  D_T(x)-D_{T'}(x) &=& \sum^{m-1}_{i=1}|V_{T_i}|\geq m-1>0\ \ \ \ \text{for all $x\in V_{G_0},$} \\
  D_T(x)-D_{T'}(x) &=& -\sum^{m-1}_{i=1}|V_{T_i}| \ \ \ \ \text{for all $x\in
V_{T_m}\cup \{v\}$}, \\
  D_T(x)-D_{T'}(x) &=& |V_{G_0}|-|V_{T_m}|-1\geq0
\ \ \ \ \text{for all $x \in V_{T_1}\cup V_{T_2}\cup \ldots\cup V_{T_{m-1}}$}. \label{eq:0036}
\end{eqnarray}

By the definition of EDS, we have
\begin{eqnarray}
  \xi^d(T)-\xi^d(T') &=& \sum_{x\in V_T}\varepsilon_T(x)D_T(x)-\sum_{x\in V_{T'}}\varepsilon_{T'}(x)D_{T'}(x) \notag\\
   &\geq & \sum_{x\in V_{T'}}\varepsilon_{T'}(x)(D_T(x)-D_{T'}(x))\ \ \ \ \text{(by (3.1) - (3.3))}\label{eq:3.1}\\
   &=& \sum_{x\in V_{G_0}}\varepsilon_{T'}(x)(D_T(x)-D_{T'}(x))+\sum_{x\in V_{T_m}\cup \{v\}}\varepsilon_{T'}(x)(D_T(x)-D_{T'}(x)) \notag\\
& & +\sum^{m-1}_{i=1}\left(\sum_{x\in V_{T_i}}\varepsilon_{T'}(x)(D_T(x)-D_{T'}(x))\right) \notag\\
   &=& \sum_{x\in V_{G_0}}\varepsilon_{T'}(x)\sum^{m-1}_{i=1}|V_{T_i}|+\sum_{x\in V_{T_m}\cup \{v\}}\varepsilon_{T'}(x)\left(-\sum^{m-1}_{i=1}|V_{T_i}|\right) \notag
\end{eqnarray}
\begin{eqnarray}
   & & +\sum^{m-1}_{i=1}\sum_{x\in V_{T_i}}\varepsilon_{T'}(x)(|V_{G_0}|-|V_{T_m}|-1)\ \ \ \text{(by (3.4)-(3.6))} \notag\\
   &\geq& \sum^{m-1}_{i=1}|V_{T_i}|\cdot\left(\sum_{x\in V_{G_0}}\varepsilon_{T'}(x)-\sum_{x\in V_{T_m}\cup \{v\}}\varepsilon_{T'}(x)\right) \ \ \ \ \ \text{(by (3.6))} \label{eq:3.2}\\
   &\geq& 0.\ \ \ \  \text{(by Fact 1)}  \label{eq:3.3}
\end{eqnarray}

The equality in (\ref{eq:3.1}) holds if and only if $\varepsilon_T(x)=\varepsilon_{T'}(x)$ for all $x\in V_T$, which is equivalent to $\varepsilon_T(v)=\varepsilon_T(w)$ (by (\ref{eq:0033})) and $|V_{P^i}| \leq |V_{T_m}|,\, i=1,2,\ldots, m-1$. 

The equality in (\ref{eq:3.2}) holds if and only if $|V_{G_0}|=|V_{T_m}|+1$,
which is equivalent to that $G_0$ is a path (otherwise, $\varepsilon_T(v)<\varepsilon_T(w)$, a contradiction).

The equality in (\ref{eq:3.3}) holds if and only if $\sum_{v\in V_{G_0}}\varepsilon_{T'}(x)=\sum_{v\in V_{T_m}\cup \{v\}}\varepsilon_{T'}(x)$.

Hence, $\xi^d(T)\ge \xi^d(T')$ with equality if and only if each of the equalities in
(\ref{eq:3.1})-(\ref{eq:3.3}) holds, i.e., $\varepsilon_T(v)=\varepsilon_T(w)$ and $T[S]$ is one of the longest paths in $T$, where
$S=V_{G_0}\cup V_{T_m}\cup \{v\}.$
\end{proof}

\begin{lem}
Suppose that $P=v_0v_1\ldots v_i\ldots v_r\ldots v_d$ is one of the longest paths contained in an $n$-vertex tree $T$ with $|V_{T_1}|\le |V_{T_{d-1}}|$ and
$r=\min\{i: \, |V_{T_i}|>1, i=2,3,\ldots, d-1\};$ see Fig. 5.  Let $T'=T-\{v_ru: u\in N_T(v_r)\setminus\{v_{r-1}, v_{r+1}\}\}+\{v_1u: u\in N_T(v_r)\setminus\{v_{r-1}, v_{r+1}\}\}.$ Then we have $\xi^d(T)<\xi^d(T').$
\end{lem}
\begin{figure}[h!]
\begin{center}
 \psfrag{t}{$v_i$} \psfrag{a}{$T_1$}\psfrag{b}{$T_2$}\psfrag{c}{$T_{r}$}\psfrag{d}{$T_{r+1}$}\psfrag{e}{$T_{d-1}$}\psfrag{m}{$T$}\psfrag{n}{$T'$}\psfrag{i}{$v_0$}
  \psfrag{j}{$v_d$}\psfrag{k}{$v_{d-1}$}\psfrag{1}{$v_1$}\psfrag{2}{$v_2$}\psfrag{r}{$v_r$}\psfrag{s}{$v_{r+1}$}
  \includegraphics[width=150mm] {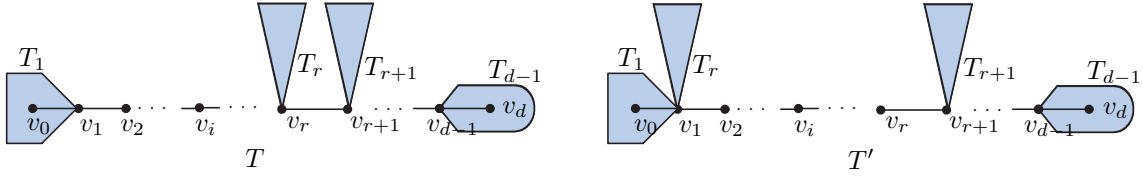}\\
  \caption{Trees $T$ and $T'$.}
\end{center}
\end{figure}
\begin{proof}
From the structure of $T$  and $T'$, it is easy to see that
$\varepsilon_T(x)\leq\varepsilon_{T'}(x)$ for $x\in V_T$.
By  the definition of EDS, we have
\begin{eqnarray*}
  \xi^d(T')-\xi^d(T)&=& \sum_{x\in V_{T'}}\varepsilon_{T'}(x)D_{T'}(x)- \sum_{x\in V_T}\varepsilon_T(x)D_T(x)\\
    &\geq & \sum_{x\in V_T}\varepsilon_T(x)(D_{T'}(x)-D_T(x)) \\
    &=& \sum_{x\in V_{T_1}}\varepsilon_T(x)(D_{T'}(x)-D_T(x)) +\sum_{i=2}^r\varepsilon_T(v_i)(D_{T'}(v_i)-D_T(v_i))\\
    & &+\sum_{x\in V_{T_r}\setminus\{v_r\}}\varepsilon_T(x)(D_{T'}(x)-D_T(x))+\sum_{x\in\bigcup_{i=r+1}^{d-1} V_{T_i}}\varepsilon_T(x)(D_{T'}(x)-D_T(x)).
\end{eqnarray*}
Note that
\begin{eqnarray*}
D_{T'}(x)-D_T(x)&=&(1-r)(|V_{T_r}|-1)\text{ for $x\in V_{T_1}$},\\
D_{T'}(v_i)-D_T(v_i)&=&(2i-r-1)(|V_{T_r}|-1),\ \ \ \text{$i=2,3,\ldots,r$},\\
D_{T'}(x)-D_T(x)&=&(1-r)|V_{T_1}|+\sum_{i=2}^r(2i-r-1)+(r-1)\sum_{i=r+1}^{d-1}|V_{T_i}| \\
&=&(r-1)(n-2|V_{T_1}|-r+3)\\
&>&0 \ \ \ \text{for $x\in V_{T_r}\setminus\{v_r\}$},\\
D_{T'}(x)-D_T(x)&=&(r-1)(|V_{T_r}|-1)\ \ \ \text{for
$x\in\bigcup_{i=r+1}^{d-1} V_{T_i}$}.
\end{eqnarray*}
Hence,
\begin{eqnarray}
 \xi^d(T')-\xi^d(T) &>&(1-r)(|V_{T_r}|-1)\sum_{x\in V_{T_1}}\varepsilon_T(x)+(|V_{T_r}|-1)\sum_{i=2}^r\varepsilon_T(v_i)(2i-r-1)\notag  \\
    && +(r-1)(|V_{T_r}|-1)\sum_{x\in\bigcup_{i=r+1}^{d-1} V_{T_i}}\varepsilon_T(x)\notag\\
    &=& (r-1)(|V_{T_r}|-1)(-\sum_{x\in V_{T_1}}\varepsilon_T(x)+\sum_{x\in\bigcup_{i=r+1}^{d-1} V_{T_i}}\varepsilon_T(x))\notag\\
    & &+(|V_{T_r}|-1)\sum_{i=2}^r\varepsilon_T(v_i)(2i-r-1).
\end{eqnarray}
It is routine to check that
$$
r-1>0,\ \ \  |V_{T_r}|-1>0,\ \ \   -\sum_{x\in
V_{T_1}}\varepsilon_T(x)+\sum_{x\in\cup_{i=r+1}^{d-1}
V_{T_i}}\varepsilon_T(x)\geq 0
$$
and
$$\sum_{i=2}^r\varepsilon_T(v_i)(2i-r-1)>\varepsilon_T(v_r)\sum_{i=2}^r(2i-r-1)=(r-1)\varepsilon_T(v_r)>0.$$
Hence, together with (3.10), we get $\xi^d(T)<\xi^d(T')$, as desired.
\end{proof}

\begin{thm}
Among $\mathscr{T}_n^k$, the balanced spider $S(\underbrace{\lceil\frac{n-1}{k}\rceil,\ldots,\lceil\frac{n-1}{k}\rceil}_r,\underbrace{\lfloor\frac{n-1}{k}\rfloor,\ldots,\lfloor\frac{n-1}{k}\rfloor}_{k-r})$
minimizes the EDS, where $n-1 \equiv r \pmod{k}$.
\end{thm}
\begin{proof}
Let $T$ be an $n$-vertex tree with $k$ leaves which has the minimal EDS, then
$T$ is a spider. Otherwise, by Lemma 3.1, there exists another $n$-vertex tree with $k$ leaves, say $T'$, such that $\xi^d(T')<\xi^d(T)$, a
contradiction. Denote $T:=S(a_1,a_2,\ldots,a_k)$. In order to complete the proof, it suffices to show that
the spider $T$ is balanced, i.e.,
that $|a_i-a_j|\leq 1, \,1\leq i,j\leq k$.

Without loss of
generality, we may suppose that $a_1\geq a_2\geq \ldots \geq
a_k\geq1$. If the spider $T$ is not balanced, then it is easy to see that $a_1-a_k \geq2$ and $\varepsilon_T(u)=a_1$, here $u$ is the hub of $T$. Denote
$P^1$ by $uu_1u_2\ldots u_{a_1-1}u_{a_1}$, while
$P^k$ by $uv_1v_2\ldots v_{a_k}$.  Let
$T_0=T-\{u_1,u_2,\ldots,u_{a_1},v_1,v_2,\ldots,v_{a_k}\}$. It is easy to see that
$\varepsilon_{T_0}(u)\le a_1$.

\textbf{Case 1}  $\varepsilon_{T_0}(u)<a_1$.

In this case, let $T'=T-\{uw: w\in N_T(u)\setminus \{u_1, v_1\}\}+\{u_1w: w\in N_T(u)\setminus \{u_1, v_1\}\}$, i.e., move $T_0$ from $u$ to $u_1$. By Lemma 3.1, we have $\xi^d(T')<\xi^d(T)$, a
contradiction.

\textbf{Case 2}  $\varepsilon_{T_0}(u)=a_1$.

Let $T'=T-u_{a_1-1}u_{a_1}+v_{a_k}u_{a_1}$ and let $\varepsilon'(x)$
(resp. $\varepsilon(x)$) be the eccentricity of $x$ in $T'$
(resp. $T$). Note that $\varepsilon(u)=a_1>a_k$ and
\begin{eqnarray*}
 && \varepsilon'(u_i)=\varepsilon(u_i)=i+a_1,  i=1,2,\ldots,a_1-1; \ \ \ \ \varepsilon'(v_j)=\varepsilon(v_j)=j+a_1,  j=1,2,\ldots,a_k,  \\
  && \varepsilon(u_{a_1})\geq \varepsilon'(u_{a_1})=1+a_1+a_k;\ \ \ \varepsilon(x)\geq\varepsilon'(x)\geq a_1,\text{ for $x\in V_{T_0}$ and
$|V_{T_0}|\geq a_1+1$}.
\end{eqnarray*}
Hence, for all $x\in V_T$, we have
\[\notag
\varepsilon(x)\geq\varepsilon'(x).
\]

Let $D'(x)$ (resp. $D(x)$) be the sum of all distances from $x$ in $T'$
(resp. $T$).
\begin{eqnarray}\label{eq:3.10}
  \xi^d(T)-\xi^d(T')&=& \Delta_1+\Delta_2+\Delta_3+\Delta_4,
\end{eqnarray}
where
$$
\begin{array}{ll}
  \Delta_1= \sum_{i=1}^{a_1-1}(\varepsilon(u_i)D(u_i)-\varepsilon'(u_i)D'(u_i)), & \Delta_2= \sum_{j=1}^{a_k}(\varepsilon(v_j)D(v_j)-\varepsilon'(v_j)D'(v_j)), \\[5pt]
  \Delta_3= \varepsilon(u_{a_1})D(u_{a_1})-\varepsilon'(u_{a_1})D'(u_{a_1}), & \Delta_4= \sum_{x\in V_{T_0}}(\varepsilon(x)D(x)-\varepsilon'(x)D'(x)).
\end{array}
$$

The contribution of vertices
$u_1,u_2,\ldots,u_{a_1-1}$ in the EDS of $T$ and $T'$ are respectively as follows
$$\sum_{i=1}^{a_1-1}(i+a_1)[(1+\cdots+(i-1))+(1+\cdots+(a_1-i))+((i+1)+\cdots+(i+a_k)+\sum_{x\in
V_{T_0}}d(x,u_i))],$$
and
$$\sum_{i=1}^{a_1-1}(i+a_1)[(1+\cdots+(i-1))+(1+\cdots+(a_1-i-1))+(i+a_k+1)+((i+1)+\cdots+(i+a_k)+\sum_{x\in
V_{T_0}}d(x,u_i))].
$$
This gives
\[\label{eq:3.11}
  \Delta_1=\frac{7}{6}a_1+\frac{3}{2}a_1a_k-a_1^2-\frac{3}{2}a_1^2a_k-\frac{1}{6}a_1^3.
\]

Similarly, the contribution of vertices
$v_1,v_2,\ldots,v_{a_k}$ in the EDS of $T$ and $T'$ are respectively as follows
$$
\sum_{j=1}^{a_k}(i+a_1)[(1+\cdots+(j-1))+(1+\cdots+(a_k-j))+((i+1)+\cdots+(j+a_1)+\sum_{x\in
V_{T_0}}d(x,v_j))]
$$
and
$$
\sum_{j=1}^{a_k}(j+a_1)[(1+\cdots+(j-1))+(1+\cdots+(a_k-j))+(a_k-j+1)+((j+1)+\cdots+(j+a_1-1)+\sum_{x\in
V_{T_0}}d(x,v_j))],
$$
which implies that
\[\label{eq:3.12}
 \Delta_2=-\frac{7}{6}a_k+\frac{1}{6}a_k^3-\frac{1}{2}a_1a_k+\frac{1}{2}a_k^2a_1+a_ka_1^2.
\]

Note that
\begin{eqnarray*}
  \varepsilon(u_{a_1})D(u_{a_1}) &=& \varepsilon(u_{a_1})\left[(1+\cdots+(a_1-1))+((a_1+1)+\cdots+(a_1+a_k))+\sum_{x\in
V_{T_0}}(d(x,u)+a_1)\right], \\
  \varepsilon'(u_{a_1})D'(u_{a_1}) &=& \varepsilon'(u_{a_1})\left[(1+\cdots+a_k)+((a_k+2)+\cdots+(a_1+a_k))+\sum_{x\in
V_{T_0}}(d(x,u)+a_k+1)\right],\\
\sum_{x\in V_{T_0}}\varepsilon(x)D(x)&=&\sum_{x\in V_{T_0}}\varepsilon(x)\left(\sum_{y\in V_{T_0}}d(y,x)+\sum_{i=1}^{a_1}(d(x,u)+i)+\sum_{j=1}^{a_k}(d(x,u)+j)\right),\\
\sum_{x\in V_{T_0}}\varepsilon'(x)D'(x)&=&\sum_{x\in V_{T_0}}\varepsilon'(x)\left(\sum_{y\in V_{T_0}}d(y,x)+\sum_{i=1}^{a_1-1}(d(x,u)+i)+\sum_{j=1}^{a_k+1}(d(x,u)+j)\right).
\end{eqnarray*}
Hence,
\begin{eqnarray}
  \Delta_3 &\geq & (|V_{T_0}|-1)(a_1+a_k+1)(a_1-a_k-1), \label{eq:3.13}\\
 \Delta_4 &\geq & \sum_{x\in V_{T_0}}\varepsilon'(x)(a_1-a_k-1)\geq a_1(a_1-a_k-1)|V_{T_0}|.\label{eq:3.14}
\end{eqnarray}

In view of (\ref{eq:3.10})-(\ref{eq:3.14}), we have
\begin{eqnarray*}
  \xi^d(T)-\xi^d(T') &\geq& \Delta_1+\Delta_2+a_1(a_1+a_k+1)(a_1-a_k-1)+a_1(a_1-a_k-1)(a_1+1)\\
   &=&(a_1-a_k-1)(\frac{11}{6}a_1^2-\frac{1}{6}a_k^2+\frac{1}{3}a_1a_k+\frac{5}{6}a_1+\frac{1}{6}a_k).
\end{eqnarray*}
Note that $a_1-a_k-1>0$ and
$\frac{11}{6}a_1^2-\frac{1}{6}a_k^2+\frac{1}{3}a_1a_k+\frac{5}{6}a_1+\frac{1}{6}a_k>\frac{11}{6}a_1^2+\frac{1}{6}a_k^2>0$.
Therefore,  $\xi^d(T)>\xi^d(T')$, a contradiction. Hence, the spider $T$
is balanced, i.e., $|a_i-a_j|\leq 1 (1\leq i,j\leq k)$, as desired.
\end{proof}

Recall that graph $P_l(a,b)$ is obtained from $P_l$ by attaching $a$ and $b$ leaves to the endvertices of $P_l$ respectively.
\begin{thm}
Let $T$ be an $n$-vertex tree with $k$ leaves, then
$
\xi^d(T)\leq\xi^d(P_{n-k}(\lfloor\frac{k}{2}\rfloor,\lceil\frac{k}{2}\rceil))
$
with equality if and only if $T\cong P_{n-k}(\lfloor\frac{k}{2}\rfloor,\lceil\frac{k}{2}\rceil).$
\end{thm}
\begin{proof}
Let $T^*$ be the $n$-vertex tree with $k$ leaves which has the maximal EDS,
then $T^*$ is of the form $P_{n-k}(a, b)$, where $a+b=k$. Otherwise, by Lemma 3.2 there exists another $n$-vertex tree with
$k$ leaves, say $\hat{T}$, such that $\xi^d(T^*)<\xi^d(\hat{T})$, a
contradiction. By Lemma 2.4, among $\{P_{n-k}(a,b):\, a+b=k,\ a,b\geq1\}$, only $
P_{n-k}(\lfloor\frac{k}{2}\rfloor,\lceil\frac{k}{2}\rceil)$ has the largest EDS.
This completes the proof.
\end{proof}

\section{\normalsize The minimal EDS of trees with a $(p,q)$-bipartition}
\setcounter{equation}{0}
Let $G$ be a connected bipartite graph with $n$ vertices. Hence its vertex set can be partitioned into
two subsets $V_1$ and $V_2$, such that each edge joins a vertex in $V_1$ with a vertex in $V_2$. Suppose that $V_1$ has
$p$ vertices and $V_2$ has q vertices, where $p+q = n$. Then we say that $G$ has a $(p, q)$-\textit{bipartition} $(p \le q)$. Let $\mathscr{T}_n^{p,q}$ be the set of all $n$-vertex trees, each of which has a $(p, q)$-bipartition $(p+q=n)$.

In this section, we are to determine the trees with the first, second and third minimal EDS in $\mathscr{T}_n^{p,q}$. Note that
$\mathscr{T}_n^{1,n-1}$ contains just $S_n$ ;
$\mathscr{T}_n^{2,n-2}=\{P_3(a, b),\, a+b=n-3\}$, where $P_3(a,b)$
is obtained from $P_3$ by attaching $a$ and $b$ leaves to the endvertices of $P_3$ respectively. By Lemmas 2.4 and 2.5, we have
$\xi^d(P_3(0,n-3))<^{\text{(Lemma 2.5)}}\xi^d(P_3(1,n-4))<\xi^d(P_3(2,n-5))<\xi^d(P_3(3,n-6))<\cdots<\xi^d(P_3(\lfloor\frac{n-3}{2}\rfloor,\lceil\frac{n-3}{2}\rceil))$.
Hence in what follows we consider $p\ge 3.$
\begin{figure}[h!]
\begin{center}
  \psfrag{a}{$w$}\psfrag{b}{$u$}\psfrag{c}{$v$}\psfrag{d}{$v_1$}\psfrag{e}{$v_2$}\psfrag{f}{$v_t$}\psfrag{u}{$T_0$}\psfrag{g}{$T$}\psfrag{h}{$T^*$}
  \includegraphics[width=120mm] {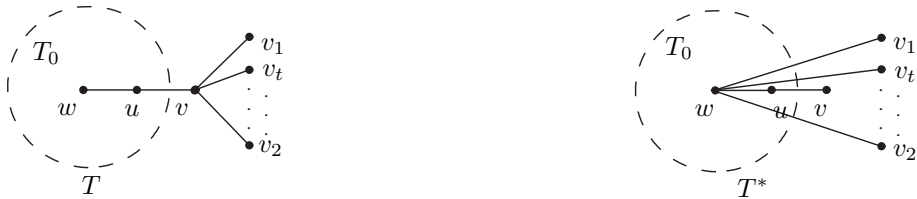}\\
  \caption{Trees $T$ and $T^*$}
\end{center}
\end{figure}
\begin{lem}
Given an $n$-vertex tree $T$ with $wu,uv\in E_T$, $d_T(w)\geq 2$, and each member
in $N(v)\setminus\{u\}=\{v_1, v_2,\ldots,v_t\}$ is a leaf, $t\geq 1$. Let $T_0=T-\{v,v_1,v_2,\ldots,v_t\}$ and
$T^*=T-\{vv_1,vv_2,\ldots,vv_t\}+\{wv_1,wv_2,\ldots,wv_t\}$. Trees $T, T_0$ and $T^*$ are depicted in Fig. 6. Then
$\xi^d(T)>\xi^d(T^*)$.
\end{lem}
\begin{proof}
It is easy to see that $\varepsilon_T(x)\geq\varepsilon_{T^*}(x)$ for
$x\in V_T$. Let $T_u$ be the component of $T-\{w,v\}$ which
contains vertex $u$. By simple calculations, we have
\begin{eqnarray*}
  D_T(x)-D_{T^*}(x)&=&0 \ \ \ \ \text{for $x\in V_{T_u}$};\ \ \ \ D_T(v_i)-D_{T^*}(v_i)=2(n-t-|V_{T_u}|-2),i=1,2,\ldots,t;\\
  D_T(v)-D_{T^*}(v)&=&-2t;\ \ \ \ \ \ \ \ \ \ \ \ \ \ \ \ \ \ D_T(x)-D_{T^*}(x)=2t\ \ \ \text{for $x\in V_{T_0}\setminus V_{T_u}$.}
\end{eqnarray*}
Therefore,
\begin{eqnarray}
 \xi^d(T)-\xi^d(T^*) &=&\sum_{x\in V_{T}}(\varepsilon_T(x)D_T(x)-\varepsilon_{T^*}(x)D_{T^*}(x))\notag\\
   &\geq& \sum_{x\in V_T}\varepsilon_{T^*}(x)(D_T(x)-D_{T^*}(x))\notag \\
   &=& 0\cdot\sum_{x\in V_{T_u}}\varepsilon_{T^*}(x)+2(n-t-|V_{T_u}|-2)\sum_{i=1}^t\varepsilon_{T^*}(v_i)-2t\varepsilon_{T^*}(v)+2t\sum_{x\in V_{T_0}\setminus V_{T_u}}\varepsilon_{T^*}(x)\notag
\end{eqnarray}
\begin{eqnarray}
      &=&2t\left((n-t-|V_{T_u}|-2)\varepsilon_{T^*}(v_1)-\varepsilon_{T^*}(v)+\sum_{x\in V_{T_0}\setminus V_{T_u}}\varepsilon_{T^*}(x)\right).\label{eq:4.1}
\end{eqnarray}
It is easy to check that $-\varepsilon_{T^*}(v)+\sum_{x\in
V_{T_0}\setminus V_{T_u}}\varepsilon_{T^*}(x)\geq 0$. Note that $d_T(w)\geq2$,
hence $n-t-|V_{T_u}|-2\geq|N_T(w)|-1\geq1>0$.
In view of (\ref{eq:4.1}), we get $\xi^d(T)-\xi^d(T^*)>0$, as desired.
\end{proof}

Note that, in Lemma 4.1, if $T$ is in $\mathscr{T}_n^{p,q}$, it is easy to see that $T^*$ is also in $\mathscr{T}_n^{p,q}$. Furthermore, ${\rm diam}(T^*)\le {\rm diam} (T)$. We call the transformation in Lemma 4.1 as
{\bf Transformation I}. Applying {\bf Transformation I} repeatedly yields the following theorem.
\begin{thm}
The tree $T(p,q)$ is the unique tree in $\mathscr{T}_n^{p,q}$ which
has the minimal EDS, where $T(p,q)$ is depicted in Fig. 7.
\end{thm}

Next, we are to determine the unique tree with the second
minimal EDS in $\mathscr{T}_n^{p,q}$. Let $\mathscr{A}=\{T_s:\, 1\leq s\leq \frac{p-1}{2}\}\bigcup \{T'_t:\, 1\leq t\leq
\frac{q-1}{2}\}$, where $T_s$ and $T'_t$ are depicted in Fig. 7.
\begin{thm}
Among $\mathscr{T}_n^{p,q}$, $T_1$ is the unique tree with the second minimal
EDS for $3\leq p\leq q$.
\end{thm}
\begin{proof}
Choose $T\in \mathscr{T}_n^{p,q}\setminus\{T(p,q)\}$ such that its EDS is as small as possible.
Note that {\bf Transformation I} strictly decreases the EDS of trees. It is easy to see that applying {\bf Transformation I} once to $T$, the resultant graph is just $T(p,q).$ Together with the definition of $\mathscr{A}$ we know the tree among $\mathscr{T}_n^{p,q}$ with the second minimal EDS must be in $\mathscr{A}$.
\begin{figure}[h!]
\begin{center}
  \psfrag{a}{$p-s-1$}\psfrag{b}{$q-2$}\psfrag{c}{$T_s$}\psfrag{d}{$s$}\psfrag{e}{$p-2$}\psfrag{f}{$q-t-1$}\psfrag{g}{$T'_t$}\psfrag{h}{$t$}
  \psfrag{r}{$p-1$} \psfrag{q}{$q-1$} \psfrag{p}{$T(p, q)$}
  \includegraphics[width=120mm]{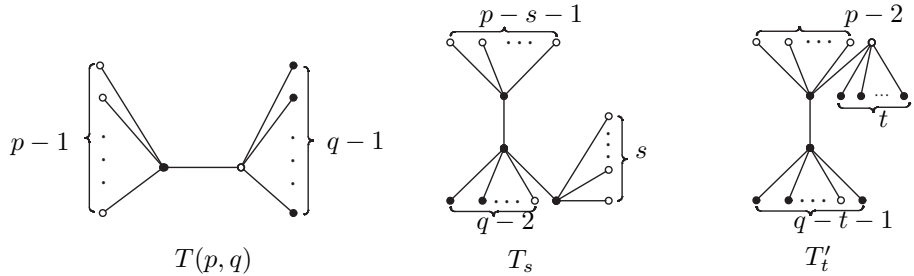}\\
  \caption{Trees $T(p,q), T_s$ and $T'_t$}
\end{center}
\end{figure}

By the definition of EDS, we have
\begin{align*}
  \xi^d(T_s) =& 4(p-s-1)(1+2(p-s-1)+3(q-1)+4s)+3(p-s+2(q-1)+3s)+2(q+2(p-1))+3(q-2)(1 \\
     &+2(q-1)+3(p-1))+3(s+1+2(q-1)+3(p-s-1))+4s(1+2s+3(q-1)+4(p-s-1))\\
   =& 6n^2+9np-7p^2-22n-4p+16ps-16s^2-16s+18=f(s).
\end{align*}
By direct verification, it follows
$f'(s)=16p-32s-16=16(p-1-2s)\geq0,$
which implies $f(s)$ is an increasing function in $s$ with $1\leq s\leq
\frac{p-1}{2}$. Hence, we have
\[\label{eq:4.2}
   \xi^d(T_1)<\xi^d(T_2)<\cdots<\xi^d(T_{\lfloor\frac{p-1}{2}\rfloor}).
\]

Similarly, we have
\begin{align*}
  \xi^d(T'_t) =& 4(q-t-1)(1+2(q-t-1)+3(p-1)+4t)+3(q-t+2(p-1)+3t)+2(p+2(q-1))+3(p-2)(1 \\
    &+2(p-1)+3(q-1))+3(t+1+2(p-1)+3(q-t-1))+4t(1+2t+3(p-1)+4(q-t-1))\\
   =& 6n^2+9nq-7q^2-22n-4q+16qt-16t^2-16t+18=g(t).
\end{align*}
By direct verification, it follows
$g'(t)=16q-32t-16=16(q-1-2t)\geq0,$
which implies $g(t)$ is an increasing function in $t$ with $1\leq t\leq
\frac{q-1}{2}$. Hence,
\[\label{eq:4.3}
   \xi^d(T'_1)<\xi^d(T'_2)<\cdots<\xi^d(T'_{\lfloor\frac{q-1}{2}\rfloor}).
\]

In order to characterize the tree with second minimal EDS
in $\mathscr{A}$, in view of (\ref{eq:4.2}) and (\ref{eq:4.3}) it suffices to compare the EDS of $T_1$ with that of $T'_1$. On the one hand, if $p=q$ we have $T_1\cong T_1'$ our result holds in this case. On the other hand, if $p<q$, by direct computing we have
\begin{eqnarray*}
   &&\xi^d(T_1)=6n^2+9np-7p^2-22n+12p-14, \ \ \ \ \xi^d(T'_1)=6n^2+9nq-7q^2-22n+12q-14.
\end{eqnarray*}
This gives that $\xi^d(T_1)-\xi^d(T'_1)=2(n+6)(p-q)<0,$ i.e., $\xi^d(T_1)<\xi^d(T'_1)$, as desired.
\end{proof}
\begin{figure}[h!]
\begin{center}
  \psfrag{a}{$p-3$}\psfrag{b}{$q-2$}\psfrag{c}{$T_2$}\psfrag{d}{$p-2$}\psfrag{e}{$q-2$}\psfrag{f}{$T'_1$}\psfrag{g}{$\hat{T}_s$}
  \psfrag{h}{$p-2$}\psfrag{j}{$q-t-2$}\psfrag{k}{$t$}\psfrag{i}{$\tilde{T}_t$}\psfrag{m}{$\vec{T}_r$}\psfrag{0}{$s$}
  \psfrag{1}{$p-s-2$}\psfrag{2}{$q-2$}\psfrag{3}{$p-3$}\psfrag{5}{$q-r-2$}\psfrag{4}{$r$}
  \includegraphics[width=160mm] {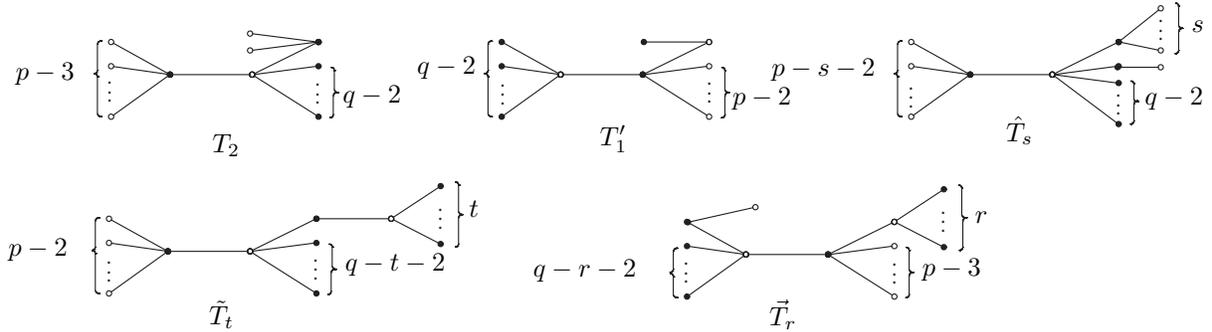}\\
  \caption{Trees $T_2, T_1', \hat{T}_s, \tilde{T}_t$ and $\vec{T}_r$.}
\end{center}
\end{figure}

Finally, we are to determine the tree with the third minimal EDS in $\mathscr{T}_n^{p,q}$.
Let $\mathscr{B}=\{T_2,\ T'_1\}\bigcup \{\hat{T}_s:\, 1\leq s\leq p-3\}\bigcup \{\tilde{T}_t:\, 1\leq t\leq p-3\}\bigcup \{\vec{T}_r:\, 1\leq r\leq
q-3\}$, where $T_2,\ T'_1, \hat{T}_s,\, \tilde{T}_t$ and $\vec{T}_r$ are depicted in Fig. 8.
\begin{thm}
Among $\mathscr{T}_n^{p,q}$ with $4\leq p< q$.
\begin{wst}
  \item[{\rm (i)}] If $n>p-3+\sqrt{p^2+9p-23}$, then $T_2$ is the unique tree with the third minimal EDS;
  \item[{\rm (ii)}] If $n<p-3+\sqrt{p^2+9p-23}$, then $T_1'$ is the unique tree with the third minimal EDS.
\end{wst}
\end{thm}
\begin{proof}
Choose $T\in \mathscr{T}_n^{p,q}\setminus\{T(p,q),\ T_1\}$ such that its EDS is as small as possible.
It is easy to see that applying {\bf Transformation I} once to $T$, the resultant graph is $T(p,q)$ or $T_1.$
On the one hand, in view of the proof of Theorem 4.3, we know $T'_1$ or $T_2$ may be the tree with the third minimal EDS among $\mathscr{T}_n^{p,q}$; on the other hand, applying {\bf Transformation I} to $T$ may yield the graph $T_1$.  Note that {\bf Transformation I} strictly decreases the EDS, hence the tree among$\mathscr{T}_n^{p,q}$ with the third minimal EDS must be in $\mathscr{B}$.

By the definition of EDS, we have
\begin{align*}
  \xi^d(T_2) =& 6n^2+9np-7p^2-22n+28p-78; \\
  \xi^d(T'_1) =& 6n^2+9nq-7q^2-22n+12q-14;\\
  \xi^d(\hat{T}_s) =& 4(p-s-2)(1+2(p-s-2)+3(q-1)+4(s+1))+3(p-s-1+2(q-1)+3(s+1)) \\
              &+2(q+2(p-1))+3(q-3)(1+2(q-1)+3(p-1))+3(2\cdot1+2(q-1)+3(p-2))+3(s+1\\
              &+2(q-1)+3(p-s-1))+4(1+2+3(q-1)+4(p-2)+4s(1+2s+3(q-1)+4(p-s-1))\\
             =&6n^2+9np-7p^2-22n+12p+16ps-16s^2-32s-14=f_1(s);\\
\xi^d(\tilde{T}_t) =& 5t(1+2t+3+4(q-t-1)+5(p-2))+4(t+1+2+3(q-t-1)+4(p-2)) \\
              &+4(q-t-2)(1+2(q-t-1)+3(p-1)+4t)+3(2\cdot1+2(q-1)+3(p-2))
\end{align*}
\begin{align*}
              &+3(q-t+2(p-1)+3t)+4(p-1+2(q-t-1)+3+4t)\\
              &+5(1+2(p-2)+3(q-t-1)+4+5t)\\
            =&8n^2+11np-9p^2-33n+14p+20nt+3pt-64t-18t^2-4=f_2(t);\\
\xi^d(\vec{T}_r) =& 5r(1+2r+3(p-2)+4(q-r-1)+5)+4(r+1+2(p-2)+3(q-r-1)+4) \\
              &+4(p-3)(1+2(p-2)+3(q-1)+4)+3(p-1+2(q-1)+3)+3(q-r+2(p-1)+3r)\\
              &+4(q-r-2)(1+2(q-r-1)+3(p-1)+4r) + 4(2\cdot1+2(q-r-1)+3(p-2)+4r)\\
              &+5(1+2+3(q-r-1)+4(p-2)+5r)\\
             =&8n^2+8np-8p^2-24n+17p+20nr-17pr-4r-18r^2-22=f_3(r).
\end{align*}

Note that $4\leq p<q,1\leq s,t\leq p-3,1\leq
r\leq q-3$, hence it follows
\begin{eqnarray}
f_1(s)&=&6n^2+9np-7p^2-22n+12p+16ps-16s^2-32s-14\geq f_1(1)\label{eq:4.4};\\
f_2(t)&=&8n^2+11np-9p^2-33n+14p+20nt+3pt-64t-18t^2-4\geq f_2(1)\label{eq:4.5};\\
f_3(r)&=&8n^2+8np-8p^2-24n+17p+20nr-17pr-4r-18r^2-22\geq f_3(1)\label{eq:4.6}.
\end{eqnarray}
In order to characterize the tree with the third minimal EDS
in $\mathscr{B}$, in view of (\ref{eq:4.4})-(\ref{eq:4.6}), it is sufficient to compare the EDS of $T'_1, T_2$, $\hat{T}_1, \tilde{T}_1$ and $\vec{T}_1$.

Note that
$$
\begin{array}{ll}
  \xi^d(\vec{T}_1)-\xi^d(\tilde{T}_1)=p^2-3np+9n-17p+42>0; & \xi^d(\hat{T}_1)-\xi^d(T_2)=16>0;\\[5pt]
  \xi^d(\tilde{T}_1)-\xi^d(\hat{T}_1)=2n^2+2np-2p^2+9n-11p-24>0. &
\end{array}
$$
Hence, it suffices for us to compare $\xi^d(T_1')$ with $\xi^d(T_2)$. In fact,
$\xi^d(T'_1)-\xi^d(T_2)= 2n^2+12n-4np-30p+64$, which yields that $\xi^d(T'_1)>\xi^d(T_2)$ if $n>p-3+\sqrt{p^2+9p-23}$ and
$\xi^d(T'_1)<\xi^d(T_2)$ otherwise.

This completes the proof.
\end{proof}

\end{document}